\documentclass[reqno,12pt]{amsart}
\usepackage{amscd,amssymb}
 
\numberwithin{equation}{section}

\theoremstyle{plain}
\newtheorem{Thm}[subsection]{Theorem}
\newtheorem{Cor}[subsection]{Corollary}
\newtheorem{Lem}[subsection]{Lemma}
\newtheorem{Prop}[subsection]{Proposition}

\theoremstyle{definition}
\newtheorem{Def}[subsection]{Definition}

\theoremstyle{remark}

\newtheorem{Rem}[subsection]{Remark}

\renewcommand{\rm}{\normalshape}

\newif\ifShowLabels
\ShowLabelstrue
\newdimen\theight
\def\TeXref#1{%
	\leavevmode\vadjust{\setbox0=\hbox{{\tt
		\quad\quad  {\small \rm #1}}}%
	\theight=\ht0
	\advance\theight by \lineskip
	\kern -\theight \vbox to
	\theight{\rightline{\rlap{\box0}}%
	\vss}%
	}}%

\ShowLabelsfalse

\renewcommand{\sec}[2]{\section{#2}\label{S:#1}%
	\ifShowLabels \TeXref{{S:#1}} \fi}
\newcommand{\ssec}[2]{\subsection{#2}\label{SS:#1}%
	\ifShowLabels \TeXref{{SS:#1}} \fi}

\newcommand{\refs}[1]{Section ~\ref{S:#1}}
\newcommand{\refss}[1]{Section ~\ref{SS:#1}}
\newcommand{\reft}[1]{Theorem ~\ref{T:#1}}

\newcommand{\refe}[1]{\eqref{E:#1}}

\newenvironment{thm}[1]%
	{ \begin{Thm} \label{T:#1}  \ifShowLabels \TeXref{T:#1} \fi }%
	{ \end{Thm} }

\renewcommand{\th}[1]{\begin{thm}{#1} \sl }
\renewcommand{\eth}{\end{thm} }

\newenvironment{lemma}[1]%
	{ \begin{Lem} \label{L:#1}  \ifShowLabels \TeXref{L:#1} \fi }%
	{ \end{Lem} }
\newcommand{\lem}[1]{\begin{lemma}{#1} \sl}
\newcommand{\elem}{\end{lemma}}

\newenvironment{propos}[1]%
	{ \begin{Prop} \label{P:#1}  \ifShowLabels \TeXref{P:#1} \fi }%
	{ \end{Prop} }
\newcommand{\prop}[1]{\begin{propos}{#1}\sl }
\newcommand{\eprop}{\end{propos}}

\newenvironment{corol}[1]%
	{ \begin{Cor} \label{C:#1}  \ifShowLabels \TeXref{C:#1} \fi }%
	{ \end{Cor} }
\newcommand{\cor}[1]{\begin{corol}{#1} \sl }
\newcommand{\ecor}{\end{corol}}

\newenvironment{defeni}[1]%
	{ \begin{Def} \label{D:#1}  \ifShowLabels \TeXref{D:#1} \fi }%
	{ \end{Def} }
\newcommand{\defe}[1]{\begin{defeni}{#1} \sl }
\newcommand{\edefe}{\end{defeni}}

\newenvironment{remark}[1]%
	{ \begin{Rem} \label{R:#1}  \ifShowLabels \TeXref{R:#1} \fi }%
	{ \end{Rem} }
\newcommand{\rem}[1]{\begin{remark}{#1}}
\newcommand{\erem}{\end{remark}}

\newcommand{\eq}[1]%
	{ \ifShowLabels \TeXref{E:#1} \fi 
	   \begin{equation} \label{E:#1} }
\newcommand{\eeq}{\end{equation}}

\newcommand{\prf}{ \begin{proof} }
\newcommand{\eprf}{ \end{proof} }

\newcommand\alp{\alpha}

\newcommand\tet{\theta}		
\newcommand\iot{\iota}

\newcommand\lam{\lambda}

\newcommand\ome{\omega}		\newcommand\Ome{\Omega}

\newcommand\calF{{\mathcal{F}}}

\newcommand\calL{{\mathcal{L}}}

\newcommand\calO{{\mathcal{O}}}

\newcommand\RR{\mathbb{R}}

\newcommand\ZZ{\mathbb{Z}}

\newcommand\CC{\mathbb{C}}

\newcommand\nek{,\ldots,}
\newcommand\sdp{\times \hskip -0.3em {\raise 0.3ex
\hbox{$\scriptscriptstyle |$}}} 

\newcommand\ind{\operatorname{ind}}

\newcommand\Ker{\operatorname{Ker}}

\newcommand\supp{\operatorname{supp}}

 \renewcommand{\d}{\text{\( \partial\)}}\newcommand{\pdrb}{\bar{\d}}

\newcommand{\F}{\calF}
\newcommand{\Fak}{\calF_a^k}
\newcommand{\n}{\nabla}
\newcommand{\hme}{{H^p(M,\calO(E))}}\newcommand{\hma}{{H^*(M,\calO(E))}}

\newcommand{\ha}{{(1,0)}}
\newcommand{\ah}{{(0,1)}}
\newcommand\ch{\operatorname{char}}
\newcommand\mult{\operatorname{mult}}

\newcommand{\lap}{\bar\square_t}

\newcommand{\ka}{K\"ahler }
\begin{document}

\title{Holomorphic Morse Inequalities and Symplectic Reduction}
\author{Maxim Braverman}
\address{Department of Mathematics\\
        Ohio State University\\
        Columbus, Ohio 43210}
\email{maxim@math.ohio-state.edu}

\begin{abstract}
We introduce Morse-type inequalities for a holomorphic circle action
on a holomorphic vector bundle over a compact \ka manifold. Our
inequalities produce bounds on the multiplicities of weights occurring
in the twisted Dolbeault cohomology in terms of the data of
the fixed points and of the symplectic reduction. This result
generalizes both Wu-Zhang extension of Witten's holomorphic Morse
inequalities and Tian-Zhang Morse-type inequalities for symplectic
reduction.

As  an application we get a new proof of the Tian-Zhang relative index
theorem for symplectic quotients.
\end{abstract}
\maketitle

\sec{introd}{Introduction}
Let $(M,\ome)$ be a compact \ka manifold of complex dimension $n$ and
let $E$ be a holomorphic Hermitian vector bundle over $M$ with the
compatible holomorphic connection.  We denote by $\hma$ the cohomology
groups with coefficients in the sheaf of holomorphic sections of $E$,
calculated from the twisted Dolbeault complex
$(\Ome^{0,*}(M,E),\pdrb_E)$.

Suppose that the circle group $S^1$ acts holomorphically and
effectively on $M$ preserving the \ka structure and that the action
can be lifted to $E$ preserving the Hermitian structure on $E$. Then,
for each $p=0,1\nek n$, we obtain a representation of $S^1$ on
$\hme$. The holomorphic inequalities we are going to discuss in this
paper give an estimate on the characters of these representations.

More precisely, assume that $\mu:M\to \RR$ is a momentum map for the
circle action on $M$ (it always exists provided the set of fixed
points is not empty \cite{Frankel}) and that $a\in\RR$ is a regular
value of $\mu$. Our inequalities estimate the character of $\hme$ in
terms of the fixed-point data and the structure of the {\em reduced
space} $M_a= \mu^{-1}(a)/S^1$.

Our inequalities contain as special cases both the Wu-Zhang extension
\cite{WuZhang} of the Witten holomorphic Morse inequalities
\cite{Witten84,MaWu} and the Tian-Zhang \cite{TianZhang1} Morse type
inequalities for geometric quantization. In fact, our inequalities may
be considered as a combination of the results of
\cite{WuZhang,TianZhang1}  (note, however, that in \cite{TianZhang1}
the inequalities are obtained for a much more general case where the
circle is replaced by an arbitrary compact Lie group). 

The proof of our main theorem is based on Witten deformation of the
Dolbeault operator $\pdrb$. Technically it is very simple since all
the necessarily calculations are already contained in
\cite{WuZhang,TianZhang1}. Moreover, we need only a very simple
version of the calculations in \cite{TianZhang1} since we work with a
circle and not with an arbitrary compact Lie group.

As an application of our inequalities we get a new proof of the
Tian-Zhang relative index theorem for symplectic quotients
\cite[Theorem~5.7]{TianZhang2} (In \cite{TianZhang2}, the result is
obtained for arbitrary symplectic manifolds. Though we prove the
theorem only for \ka manifolds, our method can be easily extended to
the case of an arbitrary symplectic manifold).


\sec{main}{Main results}

In this section we formulate our main result (\reft{main})
and discuss various applications. The proof of \reft{main} is
postponed to the next section.

\ssec{char}{Weights and formal characters} Irreducible representation
of the circle group $S^1=\{e^{i\tet}:\, \tet\in\RR\}$ are classified
by integer {\em weights} (here we use the identification of the Lie
algebra of $S^1$ with $\RR$ which takes the {\em negative} primitive
lattice element, $-2\pi i\in i\RR= Lie(S^1)$, to $1$). A representation
of weight $k\in \ZZ$ is isomorphic to the complex line $\CC$ on which
the element $e^{i\tet}\in S^1$ act by multiplications by $e^{-ik\tet}$.

If $W$ is a finite dimensional representation of $S^1$ we denote by
$\mult_k(W)$ the multiplicity of the weight $k\in\ZZ$ in $W$.
Note that the multiplicity $\mult_0(W)$ of the zero weight is equal to
the dimension of the space $W^{S^1}$ of vectors in $W$ which are
invariant with respect to the action of $S^1$.
Let {\em support} of $W$ be $\supp(W)= \{k\in\ZZ:\, \mult_k(W)\not=0\}$.

The {\em formal character} of $W$ is the formal sum 
$$
	\ch(W) \ = \ \sum_{k\in\ZZ}\mult_k(W)e^{-ik\tet}.
$$ 
It lies in the ring
$\ZZ[e^{i\tet},e^{-i\tet}]$ of Laurent polynomials in $e^{i\tet}$ with
integer coefficients. This ring is called the {\em ring of formal
characters} of the circle group. We will also consider the completion
$\calL=\ZZ[[e^{i\tet},e^{-i\tet}]]$ of this ring. The elements of
$\calL$ are formal infinite sums $q(\tet)= \sum_{k\in\ZZ}q_ke^{-ik\tet}$
where $q_k\in \ZZ$.

\ssec{mom-red}{Momentum map and symplectic reduction} Let $V$  denote
the vector field on $M$ that generates the $S^1$-action.  We will
assume that $S^1$-action is {\em Hamiltonian}, i.e. there is a moment
map $\mu\colon M\to\RR$ such that $\iot_V\ome=d\mu$. Note
(\cite{Frankel}) that it is always the case if the fixed-point set of
$S^1$ on $M$ is non-empty.

Assume that $a\in \RR$ is a regular value of the momentum map $\mu$. Then
$\mu^{-1}(a)\subset M$ is a smooth submanifold endowed with a locally
free action of $S^1$. For simplicity, we will assume that this action
is free. Then the quotient space $M_a=\mu^{-1}(a)/S^1$ is a smooth \ka
manifold. The vector bundle $E$ descends to a holomorphic Hermitian
vector bundle $E_a$ over $M_a$. 

Let $\calF=M\times\CC$ denote the trivial line bundle over $M$ with
$S^1$ action defined by the formula 
\eq{F}
	e^{i\tet}:\, (x,z) \ \mapsto \ (e^{i\tet}\cdot x, e^{i\tet}z),
	\qquad x\in M, z\in \CC.
\end{equation}
Denote by $q:\mu^{-1}(a)\to M_a=\mu^{-1}(a)/S^1$ the projection. The
restriction $\calF|_{\mu^{-1}(a)}$ of $\calF$ on $\mu^{-1}(a)$
descends to a unique bundle $\calF_a$ over $M_a$ such that
$q^*\calF_a=\calF|_{\mu^{-1}(a)}$.

Let \/ $\calF^{-1}, \calF_a^{-1}$ \/ be the dual bundles to $\calF$
and \/ $\calF_a$ \/ respectively. If \/ $k\ge 0$, we denote by \/
$\F^{\pm k}$ \/ (resp.  \/ $\F_a^{\pm k}$) the $k$-th tensor power of the
bundle \/ $\F^{\pm1}$ \/ (resp. \/ $\F_a^{\pm1}$).  Note that \/
$q^*\Fak= \F^k$. Obviously,
\eq{multF}
	\mult_m H^{*}(M,\calO(E\otimes \F^k)) \ = \
		\mult_{m+k} H^{*}(M,\calO(E))
\end{equation}
for any $k,m\in\ZZ$.

We will be interested in cohomology \/
$H^{*}(M_a,\calO(E_a\otimes \Fak))$ of $M_a$ with coefficients in the
sheaf of holomorphic sections of $E_a\otimes \Fak$.

There is another action of $S^1$ on $\F$ given by the formula
\eq{action'}
	e^{i\tet}:\, (x,z) \ \mapsto \ (x, e^{i\tet}z),
	\qquad x\in M, z\in \CC.
\end{equation}
This action commutes with \refe{F} and, hence, reduces to
$\F_a$. Therefore we obtain induced actions of $S^1$ on $\F_a^k$. This
action preserve the base points in $F$ and the weight of the
representation on the fibers is equal to \/ $-k$.

Let $S(\F_a,\F_a^{-1})$ denote the direct sum 
$$
	S(\F_a,\F_a^{-1}) \ = \ \bigoplus_{k\in\ZZ} \F_a^k.
$$
One should think about \/ $S(\F_a,\F_a^{-1})$ \/ as about symmetric algebra
in $\F_a,\F_a^{-1}$.
This is an infinite-dimensional $S^1$-equivariant bundle over
$M_a$. However the cohomology groups \/ $H^*(M_a,\calO(E_a\otimes
S(\F_a,\F_a^{-1})))$ \/ are sums of those with coefficients in
$\F_a^k$. Therefore the multiplicity of any weight $k$ in \/
$H^*(M_a,\calO(E_a\otimes S(\F_a,\F_a^{-1})))$ \/ is finite and the
character 
\eq{chS}
	\ch H^*(M_a,\calO(E_a\otimes S(\F_a,\F_a^{-1}))) \ = \
		\sum_{k\in\ZZ}e^{ik\tet} 
		\dim_\CC H^{*}(M_a,\calO(E_a\otimes \Fak)).
\end{equation}
is well defined as an element of $\calL$.

\ssec{fixed}{Fixed-point set} Suppose that $F$ is a connected
component of the fixed-point set of the $S^1$-action on $M$. Then $F$
is a compact \ka manifold.  Let  $n-n_F$ be  the complex dimension
of $F$. The complexification of the normal bundle $N_F\to F$ in $M$
has the decomposition $N^\CC_F=N^\ha_F\oplus N^\ah_F$, where $N^\ha_F$ is
a holomorphic vector bundle over $F$ of rank $n_F$. The circle
$S^1$ acts on $N_F$ preserving the base points in $F$.  Moreover, the
weights of the isotropy representations on the normal fiber are
constant along $F$. 

Let \/ $\lam_{k}$ ($1\le k\le n-n_F$) \/ be the isotropy weights on
$N_F^\ha$. Let \/ $N^{\pm,\ha}_F$ \/ be the direct sum of the
sub-bundles corresponding to the weights $\lam_k>0$ and $\lam_k<0$
respectively. We denote by $\nu_F$ the rank of the holomorphic vector
bundle $N^{-,(1,0)}_F$.

The {\em polarized symmetric tensor products}
(cf. \cite{SilGuil96,WuZhang}) are the vector bundles
\eq{PSTP}\begin{aligned}
	K_F^+ \ &= \ S((N^{+,\ha}_F)^*) \ \otimes \ 
		S(N^{-,\ha}_F)\otimes\det(N^{-,\ha}_F), \\
	K_F^- \ &= \ S((N^{-,\ha}_F)^*) \ \otimes \ 
		S(N^{+,\ha}_F)\otimes\det(N^{+,\ha}_F).
	\end{aligned}
\end{equation}
Here \/ $S((N^{\pm,\ha}_F)^*)$, $S(N^{\pm,\ha}_F)$ \/ denote the sums of all
symmetric powers of the bundles \/ $(N^{\pm,\ha}_F)^*$ \/ and
$N^{\pm,\ha}_F$ respectively and \/ $\det(N^{\pm,\ha}_F)$ \/ denotes the top
exterior power of \/ $N^{\pm,\ha}_F$.

The fiber $E_p$ over each fixed point $p\in F$ is a representation of
$S^1$, and $\ch(E_p)$ is independent on $p\in F$. Consider infinite
dimensional holomorphic bundles $K_F^\pm\otimes E|_F$. The circle acts
on the total space while preserving the base points in $F$. A
sub-bundle of any given weight is a holomorphic vector bundle of
finite rank, i.e., \eq{KF=} K_F^\pm\otimes E|_F \ = \ \oplus_{k\in
\ZZ}E_{F,k}^\pm,
\end{equation}
where $E_{F,k}^\pm$ is a $S^1$-invariant sub-bundle of finite rank on
which the circle acts with weight $k$.
The cohomology groups $H^*(F,\calO(K_F^\pm\otimes E|_F))$ are the sum of
those with coefficients in $E_{F,k}^\pm$, each equipped with an induced
$S^1$-action. Therefore, for any $k\in\ZZ$ the multiplicities
\eq{mult}
	\mult_k H^*(F,\calO(K_F^\pm\otimes E|_F)) \ = \ \dim_\CC
	H^*(F,\calO(E_{F,k}^\pm))
\end{equation}
are finite. 

Our main result is the following
\th{main} For any $k\in \ZZ$ there exists a polynomial $Q_k(t)$ with
  non-negative integer coefficients such that 
  \begin{multline}\label{E:HMI-mult}
	\sum_{p=0}^{n-1} t^p\dim_{\CC} 
		H^{p}(M_a,\calO(E_a\otimes \Fak))  
	\ + \ 
	\sum_{\mu|_F>a} t^{n_F-\nu_F} \sum_{p=0}^{n_F} t^p\mult_k 
	H^p(F,\calO(K_F^-\otimes E|_F)) \\ \ + \ 
	\sum_{\mu|_F<a} t^{\nu_F} \sum_{p=0}^{n_F} t^p\mult_k
	H^p(F,\calO(K_F^+\otimes E|_F)) \\ 
	\ = \
	\sum_{p=0}^{n} t^p\mult_k H^p(M,\calO(E)) \ + \ (1+t)Q_k(t) 
   \end{multline}
   where the second and third sums in the left hand side are taken over
   connected components of the fixed-point set.
\eth
\reft{main} is proven in \refs{proof}.

The rest of this section is devoted to discussion of different
applications and reformulation of \reft{main}.

\ssec{char'}{Reformulation in the language of characters} 
 It follows from \refe{mult} that the character 
\eq{charK}
	\ch H^p(F,\calO(K_F^\pm\otimes E|_F)) 
	\ = \ 
	\sum_{k\in\ZZ} e^{-ik\tet}\,
	\dim_\CC H^p(F,\calO(E^\pm_{F,k}))
\end{equation}
of the infinite dimensional representation $H^p(F,\calO(K_F^\pm\otimes
E|_F))$ is well defined as an element of $\calL$.

\defe{polyn} Let  $q(\tet)= \sum_{k\in\ZZ}q_ke^{-ik\tet}$ be a formal
  character of $S^1$, we say $q(\tet)\ge 0$ if $q_k\ge 0$ for all
  $k\in\ZZ$. Let $Q(\tet,t)= \sum_{m=0}^n q_m(\tet)t^m$ be a
  polynomial of degree $n$ with coefficients in $\calL$, we say
  $Q(\tet,t)\ge0$ if $q_m(\tet)\ge0$ for all $m$.

  For two such polynomials $P(\tet,t)$ and $q(\tet,t)$, we say
  $P(\tet,t)\le Q(\tet,t)$ if the exists a polynomial
  $Q(\tet,t)-P(\tet,t)\ge0$.
\edefe

Using \refe{chS} we can reformulate \reft{main} in the language of
characters. 
\th{main'}
  There exists a polynomial $Q(\tet,t)\in \calL[t]$,
  such that $Q\ge0$ and 
  \begin{multline}\label{E:holMineq} 
	\sum_{p=0}^{n-1} t^p
		\ch H^{p}(M_a,\calO(E_a\otimes S(\F_a,\F_a^{-1})))  
	\ + \ 
	\sum_{\mu|_F>a} t^{n_F-\nu_F} \sum_{p=0}^{n_F} t^p\ch
	H^p(F,\calO(K_F^-\otimes E|_F)) \\
	\ + \ 
	\sum_{\mu|_F<a} t^{\nu_F} \sum_{p=0}^{n_F} t^p\ch
	H^p(F,\calO(K_F^+\otimes E|_F)) \\ 
	\ = \
	\sum_{p=0}^{n} t^p\ch H^p(M,\calO(E)) \ + \ (1+t)Q(\tet,t) 
   \end{multline}
\eth

\ssec{wit-zh-wu}{Witten-Wu-Zhang inequalities} \reft{main'}
provides estimates on the character of $\hma$ for any regular value
$a\in\RR$ of the momentum map. Let as choose $a< \min\{\mu(x):\, x\in
M\}$. Then the reduced space $M_a$ is empty and the first and the third
summands in the left hand side of \refe{holMineq} vanish. Hence,
\refe{holMineq} reduces to
\eq{ZW}
	\sum_{F} t^{n_F-\nu_F} \sum_{p=0}^{n_F} t^p\ch
	H^p(F,\calO(K_F^-\otimes E|_F))
	\ = \
	\sum_{p=0}^{n} t^p\ch H^p(M,\calO(E)) \ + \ (1+t)Q(\tet,t) 
\end{equation}
where the sum in the left is taken over all connected components of
the fixed-point set. This is precisely the Wu-Zhang extension of the
Witten holomorphic Morse inequalities for a circle action 
\cite[Theorem~2.4]{WuZhang}.

Note that choosing $a> \max\{\mu(x):\, x\in M\}$ leads to inequalities
which are similar but different from \refe{ZW}. It is shown in
\cite{Witten84} that combination of those inequalities with \refe{ZW}
gives much better estimates than \refe{ZW} alone. Even more
information about $\hma$ may be obtained by considering
\refe{holMineq} with all possible values of $a$.

\ssec{TZ}{The Tian-Zhang inequalities for symplectic reduction} Let
$F$ be a connected component of the fixed-pont set. As we have already
mentioned in \refss{fixed}, the circle acts on the normal bundle $N_F$
to $F$ preserving the base points in $F$. Also the weights of the
isotropy representation on the normal fiber are constant along $F$.

Let $\lam^+_F$ (resp. $\lam^-_F$) denote the sum of the positive
(resp. negative) weights. One easily checks that \/ \/ $\supp
K_F^+\subset(-\infty,-|\lam^-_F|]$ \/ and \/ $\supp
K_F^-\subset[\lam^+_F,\infty)$. It follows that, if $\supp
E|_F\subset[k_1,k_2]$ \/ then
\eq{supp}\begin{aligned}
	\supp \dim H^p(F,\calO(K_F^+\otimes E|_F)) \ 
		\subset \ (-\infty,k_2-|\lam^-_F|]; \\
 	\supp \dim H^p(F,\calO(K_F^-\otimes E|_F)) \ 
 		\subset \ [k_1+\lam^+_F,\infty).
\end{aligned}\end{equation}
If there exists an integer $k$ which is greater than \/
$k_2-|\lam^-_F|$ \/ for any $F$ with \/ $\mu(F)<a$ \/ and which is
smaller than \/ $k_1+\lam^+_F$ \/ for any $F$ with \/ $\mu(F)>a$ \/
then \refe{HMI-mult}, \refe{supp} imply that (cf. \cite[\S4]{SiKaTo})
\eq{TZ'} 
	\sum_{p=0}^{n-1} t^p\dim_{\CC} 
		H^{p}(M_a,\calO(E_a\otimes\Fak)) 
	\ = \ 
	\sum_{p=0}^{n} t^p\mult_k
		H^p(M,\calO(E)) \ + \ (1+t)Q_k(t).
\end{equation}

Consider the case when $E$ is a {\em pre-quantum line bundle}. That
means that the K\"ahler form $\ome$ represents the Chern class of $E$
in the cohomology $H^2(M)$.  For this case there is a natural choice
of the momentum map $\mu$ given by the Kostant formula: the
infinitesimal generator of the action of $S^1$ on the space of
sections of $E$ is given by $-\n^E_V+2\pi{i}\mu$, where $\n^E_V$
is the covariant derivative along the vector field $V$ which generates
the action of $S^1$ on $M$.

With this choice of $\mu$ one easily checks that the weight of the
representation of $S^1$ on $E|_F$ is positive (resp. negative) if
$\mu(F)>0$ (resp. $\mu(F)<0$). It follows that \refe{TZ'} holds for
$k=0, a=0$. In other words, there exists a polynomial $Q(t)$ with
non-negative coefficients such that
\eq{TZ}
	\sum_{p=0}^{n-1} t^p\dim_{\CC} H^{p}(M_0,\calO(E_0))  
	\ = \
	\sum_{p=0}^{n} t^p \dim H^p(M,\calO(E))^{S^1} \ + \ (1+t)Q(t).
\end{equation}
(here $H^p(M,\calO(E))^{S^1}$ denotes the space of $S^1$ invariant
vectors in $H^p(M,\calO(E))$). Equation \refe{TZ} is precisely the
Tian-Zhang Morse-type inequalities for symplectic reduction on a
K\"ahler manifold \cite[Theorem~5.1]{TianZhang1} (note, however, that
in \cite{TianZhang1} the inequalities are obtained for a much more
general case where the circle is replaced by an arbitrary compact Lie
group).

\ssec{index}{Index theorem} An interesting corollary of \reft{main}
may be obtained by setting $t=-1$ in \refe{HMI-mult}. Then the last
summand in the right hand side of \refe{HMI-mult} vanishes and we
obtain a combination of the Atiyah-Bott fixed point
theorem\cite{AtBott-L1,AtBott-L2} and the Guillemin-Sternberg
``quantization commutes with reduction'' theorem
\cite[Theorem~5.2]{GuiSter82}. We will now explain this in more
details.

For a connected component $F$ of the fixed-point set, define
\eq{ind}\begin{aligned}
	\ind_k(F;K_F^+\otimes E|_F) \ = \ 
		\sum_{p=0}^{n_F} t^{p+\nu_F}\mult_k
			H^p(F,\calO(K_F^+\otimes E|_F)); \\
	\ind_k(F;K_F^-\otimes E|_F) \ = \ 
		\sum_{p=0}^{n_F} t^{p+n_F-\nu_F}\mult_k
			H^p(F,\calO(K_F^-\otimes E|_F)).
\end{aligned}
\end{equation}
Recall that  bundles $E_{F,k}^\pm$ are introduced in \refe{KF=}.
By the Riemann-Roch-Hirzebruch theorem \cite{AtSinger1} (see also
\cite[Theorem~4.9]{BeGeVe}) we have
\eq{ind'}
	\ind_k(F;K_F^\pm\otimes E|_F) \ = \
	  \int_F\, Td(F)\, ch (E_{F,k}^\pm),
\end{equation}
where $Td$ and $ch$ stand for the Todd class and Chern character
respectively.

Setting $t=-1$ in \refe{HMI-mult} and taking into account \refe{ind},
we obtain
\begin{multline}\label{E:indth}
	\sum_{p=0}^{n-1} (-1)^p\dim_{\CC} 
		H^{p}(M_a,\calO(E_a\otimes\Fak))  \ + \ 
	\sum_{\mu|_F>a} \ind_k(F;K_F^-\otimes E|_F) 
	\\ \ + \ 
	\sum_{\mu|_F<a} \ind_k(F;K_F^+\otimes E|_F)
	\ = \
	\sum_{p=0}^{n} (-1)^p\mult_k H^p(M,\calO(E)). 
\end{multline}
If in \refe{indth} we choose \/ $a<\min\{\mu(x):\, x\in M\}$, then the
first term in the left hand side of \refe{indth} vanishes and, in
view of \refe{ind'}, we get the Atiyah-Bott fixed-point theorem. From
the other side, if $E$ is a pre-quantum line bundle $a=0$  and $k=0$, then
(cf. \refss{TZ}) the second and the third terms in the left hand side
of \refe{indth} vanish and we obtain the Guillemin-Sternberg
``quantization commutes with reduction'' theorem
\cite[Theorem~5.2]{GuiSter82}:
$$
	\sum_{p=0}^{n-1} (-1)^p\dim_{\CC} 
		H^{p}(M_0,\calO(E_0)) 
	\ = \ 
	\sum_{p=0}^{n} (-1)^p\dim_\CC H^p(M,\calO(E))^{S^1}. 
$$

\ssec{TZindex}{The Tian-Zhang relative index theorem for symplectic
quotients} Let $a<b$ be two regular values of the momentum map $\mu$
and suppose that $S^1$ acts freely on $\mu^{-1}(a), \mu^{-1}(b)$. Then
we can form smooth manifolds $M_a$ and $M_b$ as in \refss{mom-red} and
then apply \refe{indth} to each one of them. Comparing the results we
obtain
\begin{multline}\label{E:TZindex}
	\sum_{p=0}^{n-1} (-1)^p\dim_{\CC}
		H^{p}(M_b,\calO(E_b\otimes\F_b^k)) 
	\ - \
	\sum_{p=0}^{n-1} (-1)^p\dim_{\CC} 
		H^{p}(M_a,\calO(E_a\otimes\Fak)) 
	\\ = \ 
	\sum_{a<\mu|_F<b} \ind_k(F;K_F^-\otimes E|_F) \ - \
	\sum_{a<\mu|_F<b} \ind_k(F;K_F^+\otimes E|_F).
\end{multline}
This formula was first obtained by Tian and Zhang
\cite[Theorem~5.7]{TianZhang2} using rather sophisticated study of the
spectral flow of a family of Dirac operators with the
Atiyah-Patodi-Singer boundary conditions on a simplectic manifold 
with boundary. (In fact, Tian and Zhang considered only the case
$k=0$. However, \refe{TZindex} follows easily from their result).

The result of Tian and Zhang is valid for a more
general case where $M$ is an arbitrary symplectic manifold. Note that
our proof may be easily extended to that case.

\sec{proof}{Proof of \reft{main}}

First of all, note that it is enough to prove \reft{main} for $k=0$.
Indeed, suppose that the theorem is proven for $k=0$ and recall that
bundles $\F^k$ are defined in \refss{mom-red}.  Applying \reft{main}
with $k=0$ to the tensor product $E\otimes \F^m$ and using
\refe{multF} we obtain the statement of the theorem for $k=m$.

Let us prove \reft{main} for $k=0$. In other words, we will be
interested in $S^1$-invariant elements of $H^{*}(M,\calO(E))$ and 
$H^*(F,\calO(K_F^\pm\otimes E|_F))$. Also, without loss of generality,
we assume that $a=0$.

Recall that $\mu:M\to \RR$ is a momentum map for the circle action on
$M$.  Following \cite{TianZhang1}, we consider a one parameter family
of differentials $\pdrb_t:\Ome^{0,*}(M,E)\to \Ome^{0,*+1}(M,E)$
defined by
$$
	\pdrb_t\alp \ = \ e^{-t|\mu|^2}\, \pdrb\,  e^{t|\mu|^2}\alp 
	\ = \
	\pdrb\alp \ + \ 2\mu \pdrb\mu\wedge\alp.
$$
Let  $\pdrb^*_t$ denote the formal adjoint to $\pdrb_t$ and consider the
corresponding Laplacian
$$
	\lap \ = \ \pdrb_t^*\pdrb_t \ + \ \pdrb_t\pdrb_t^*.
$$
Clearly, for each $t\in\RR$ the cohomology $\hma$ is isomorphic to the
kernel $\Ker\lap$ of $\lap$. Moreover, the $S^1$ invariant part of
$\hma$ is isomorphic to the kernel of the restriction of $\lap$ on the
space $(\Ome^{0,*}(M,E))^{S^1}$ of $S^1$-invariant anti-holomorphic
differential forms. The later operator is calculated in
\cite{TianZhang1}. It is shown in \cite{TianZhang1} that, for 
$t\to\infty$, the calculation of the kernel may be localized to small
neighborhoods of $\mu^{-1}(0)$ and of fixed-point set of the action of
$S^1$. Such a localization, by standard techniques of
\cite{Witten82,Witten84,TianZhang1,WuZhang} leads to  Morse-type inequalities. 
The contribution of the $\mu^{-1}(0)$ to these inequalities is
calculated in \cite{TianZhang1} and is precisely equal to the first
summand in the left hand side of \refe{HMI-mult}.

The contribution of the fixed-point set to the inequalities may be
calculated using the technique of \cite{WuZhang}.  Let $F$ be a
connected component of the fixed-point set. Then the restriction of
$\mu$ on $F$ is a constant. Moreover, $\mu(F)\not=0$ since $0$ is a
regular value of $\mu$. 

In \cite{WuZhang}, Zhang and Wu considered a
one parameter deformation of $\pdrb$ given by 
$$
	\pdrb'_s\alp \ =\ \pdrb\alp \ + \   s\pdrb\mu\wedge\alp.
$$
Near $F$ our operator $\pdrb_t$ looks like $\pdrb_s'$ with
$\mu(F)s=t$. Hence, if $\mu(F)>0$, the asymptotic behavior for
$t\to\infty$ of the eigenforms of $\lap$ which concentrate near $F$ is
the same as in \cite{WuZhang}. In particular, the contribution of $F$
to the inequalities is the same as in \cite{WuZhang}. This leads to
the second summand in the left hand side of \refe{HMI-mult}.

If $\mu(F)<0$, then, as $t\to \infty$, the operator $\lap$ behaves 
as the corresponding operator in \cite{WuZhang} behaves for $s\to
-\infty$. This leads to the last term in the left hand side of
\refe{HMI-mult}.

\providecommand{\bysame}{\leavevmode\hbox to3em{\hrulefill}\thinspace}

\end{document}